\documentclass[11pt,reqno]{amsart}
\usepackage{amssymb,latexsym}
\usepackage{graphicx}
\usepackage{url}		%does nice formatting of URLs

\calclayout

\makeatletter
\newcommand{\monthyear}[1]{%
  \def\@monthyear{\uppercase{#1}}}
\newcommand{\volnumber}[1]{%
  \def\@volnumber{\uppercase{#1}}}
\AtBeginDocument{%
\def\ps@plain{\ps@empty
  \def\@oddfoot{\@monthyear \hfil \thepage}%
  \def\@evenfoot{\thepage \hfil \@volnumber}}
\def\ps@firstpage{\ps@plain}
\def\ps@headings{\ps@empty
  \def\@evenhead{%
    \setTrue{runhead}%
    \def\thanks{\protect\thanks@warning}%
    \uppercase{The Fibonacci Quarterly}\hfil}%
  \def\@oddhead{%
    \setTrue{runhead}%
    \def\thanks{\protect\thanks@warning}%
    \hfill\uppercase{On Higher Order K-bonacci Matrices}}%
  \let\@mkboth\markboth
  \def\@evenfoot{%
    \thepage \hfil \@volnumber}%
  \def\@oddfoot{%
    \@monthyear \hfil \thepage}%
  }%
\footskip=25pt
\pagestyle{headings}%
}
\makeatother

\theoremstyle{plain}
\numberwithin{equation}{section}
\newtheorem{thm}{Theorem}[section]
\newtheorem{theorem}[thm]{Theorem}
\newtheorem{lemma}[thm]{Lemma}

\begin{document}
%% replace the values in the next three lines by the correct information
\monthyear{Month year}
\volnumber{Volume, Number}
\setcounter{page}{1}

\title{On Higher Order K-bonacci Matrices}
\author{Shubhra Gupta}
\address{School of Business Administration\\
University and Institute of Advanced Research\\
Koba Institutional Area\\
Gandhinagar - 382426 Gujarat, India}
\email{shubhra.gupta@iar.res.in}
%\thanks{Research supported in part by the Natural Sciences and %Engineering Research Council of Canada and by Emperor Frederick %II of Sicily.}
%\author{Leonardo Pisano}
%\address{Dipartimento di Matematica\\
               %Universit\`{a} di Pisa\\
              % 144 Via Fibonacci\\
              % 56127 Pisa, Italy}
%\email{leo@dm.unipi.it}

\begin{abstract}
In this paper, we have constructed the higher order k-bonacci matrices and studied some of their basic properties. We have also shown that these matrices satisfying some new and interesting relations in k-bonacci recurrence. This is the interesting generalization of the work of Z. Cheng-Heng \cite{che01,chen01}.
\end{abstract}

\maketitle
%%%%%%%%%%%%%%%%%%%%%%%%%%%%%%%%%%%%%%%%%%%%%%%%%%%%%%%%%%%%%%%%%%%%%%%%%%%%%%%%%%%%%%%%%%%%%%%
\section{Introduction}
Fibonacci numbers have many applications, including appearance in the nature, such as the branching of trees, the patterns on a pineapple, the florets of a sunflower, the spirals of a pine cone, the family trees of cows and bees, and the placement of leaves on the stems of many plants, etc. These numbers were known in ancient India as Hemachandra numbers. The key property of the Fibonacci numbers is that every number equals to the sum of the two preceding numbers, therefore we obtain a numerical sequence ($f_0,f_1,\ldots,f_n, f_{n+1}\ldots$ such that $f_{n+1}=f_{n}+f_{n-1}$ with $f_0=1,f_1=1$) of numbers. The sequence is called the recurrent sequence and the relation is known as the recurrence relation \cite{nnm02}. One can generalize the idea to higher order recurrence such as Tribonacci sequence, where every number equals to the sum of the previous three numbers and with the sum of previous four numbers one can make tetrabonacci sequence and so forth. Similarly, in a more general settings, for the K-bonacci sequence every number equals to the sum of the previous k numbers \cite{Mil60}.

The idea of the sequence can be generalized to matrices. Many researchers have studied such interesting matrices. Fibonacci numbers have been generalized to obtain a square Fibonacci matrix [6]. Further, it is generalized to Fibonacci Q-matrix [5] and sequence of matrices of order $2^r$ by Cheng-Heng [1, 2].
In this paper, we obtain a generalization of the matrices of order $2^r$ to the k-bonacci matrices of order $k^r$ and obtain some interesting results. k-bonacci numbers and its matrices have various applications in tournament sequences, graph theory, modular form, dynamic interpretation and memory allocation schemes \cite{Lee00,Lee97}.

The paper is organized as follows. Section $2$ contains some preliminaries of the $k$-bonacci numbers. Section $3$ defines the $k$-bonacci matrices of order $k^r$. In section $4,$ we obtain some of its interesting properties. Section $5$ contains some new results of Fibonacci matrices of order $2^r$.
%%%%%%%%%%%%%%%%%%%%%%%%%%%%%%%%%%%%%%%%%%%%%%%%%%%%%%%%%%%%%%%%%%%%%%%%%%%%%%%%%%%%%%%%%%%
\section{Preliminaries and Notations}

For an integer $k \geq 2$, the $k$-bonacci numbers $f_{j,k}$ are defined as \cite{Mil60}
\begin{equation}\label{one-e}
f_{0,k}=0, f_{1,k}=0, f_{2,k}=0 \;\; ...\;\; f_{k-2,k}=0, f_{k-1,k}=1 \quad\text{and} \quad f_{j,k}=\sum_{n=1}^{k} f_{j-n,k} \quad\text{for all j} \geq k. 
\end{equation}
\noindent It is easy to obtain some values for $f_{j,k}$, as given below:
\[
\begin{array}{ccc}
f_{j,k}&=&2^{j-k},\; k \leq j \leq 2k-1\\   
       &=&2^{j-k}-(2^{j-2k+1}-1),\; 2k \leq j \leq 3k-2.       
\end{array}
\]
 
\noindent When $k=2$, the equation $(2.1)$ reduces to the usual Fibonacci numbers $f_{j,2}$ or $f_j$ (for simple notation),
\[
f_{0}=0\quad\text{for}\quad j = 0, f_{1}=1\quad\text{for}\quad j = 1\quad\text{and}\quad f_{j}=f_{j-1}+f_{j-2}\quad\text{for all j},\quad j \geq 2.
\]

\noindent The backward k-bonacci numbers are given by the following recurrence relation
\begin{equation}\label{b-k-on}
f_{j-k,k}=f_{j,k}- \sum_{n=1}^{k-1} f_{j-n,k}
\quad\text{for all j},\quad j \leq k-1.
\end{equation}

\noindent From backward recurrence relation, one can easily obtain some $k$-bonacci backward values:
\[
\begin{array}{ccc}
f_{j,k}&=&1,-1,0...0,\; 0 \leq j \leq k-1\\   
       &=&2,-3,1,0...0,\; -k \leq j \leq -1\\
       &=&4,-8,5,-1,0...0,\; -(2k+1) \leq j \leq -(k+1)\\
       &=&8,-20,18,-7,1,0...0,\; -(3k+2) \leq j \leq -(2k+2)\\
       &=&16,-48,56,-32,9,-1,0...0,\; -(4k+3) \leq j \leq -(3k+3)       
\end{array}
\]

From above values, we observe the following interesting properties in every block of $k$ terms:
\begin{enumerate}
\item Numbers are alternative positive and negative.
\item Sum of the numbers are zero.
\item First number belongs to k-bonacci numbers, for $j=nk\quad\text{where}\quad n \geq 0$. From forward and backward k-bonacci numbers, one can see that
\[
\begin{array}{ccccc}
  f_{k,k}&=&f_{-1}&= &1\\
  f_{k+1,k}&=&f_{-(k+1)}&=& 2\\
  f_{k+2,k}&=&f_{-(2k+1)}&= &4\\
 \ldots&=&\ldots&=&\ldots\\
 \end{array}
 \]
\noindent Therefore, we find an interesting relation.
\[
f_{k+j,k}=f_{-(jk+1),k}=2^{j}, j \geq 0.
\]

\item The last non zero number is $-1$ or $1$ for $j=n(k+1)+1 \quad\text{where} \quad  n \geq 0$.
\item The second last non zero number is subsequent odd number for $j=n(k+1)\quad\text{where}\quad n \geq 0$. 
Every number before that is even except for the first block (there are only two terms).
\end{enumerate}

\noindent For $k=2$, equation $(2.2)$ becomes the backward Fibonacci relation \cite{Vaj89} 
\[
f_{j-2}=f_{j}-f_{j-1}\quad\text{for all j},\quad j>1.
\]

Now we focus our attention to the matrices. The k-bonacci matrix \cite{Mil60} of order $k$ is defined as

\begin{equation}\label{3}
F_{j,k}^{(k)}=\left(\begin{array}{ccccc} 
f_{j+k-1,k} & \ldots & f_{j,k}\\f_{j+k-2,k} & \ldots & f_{j-1,k}\\ \vdots & \vdots & \vdots\\f_{j,k} & \ldots & f_{j-k+1,k} \end{array}\right)
=\left(f_{j+k-\lambda-\mu+1,k}\right)=f_{\lambda\mu}, \quad1 \leq \lambda,\mu \leq k.
\end{equation}
\noindent where $j \geq 0$.

\noindent By taking $k=2$, the above matrix $F_{j,k}^{(k)}$ reduces to Fibonacci matrix of order $2$, which is considered in \cite{che01}.

\noindent In view of equation (\ref{b-k-on}) we construct the backward k-bonacci matrix of order $k$ as
\[
F_{-j,k}^{(k)}=\left(f_{-(j+k-\lambda-\mu+1),k}\right), j \geq 0.
\]

%%%%%%%%%%%%%%%%%%%%%%%%%%%%%%%%%%%%%%%%%%%%%%%%%%%%%%%%%%%%%%%%%%%%%%%%%%%%%%%%%%%%%%%%%%%%%%%
\section{Matrices of Order $k^{r}$}

We have constructed the $k^{r}$th order $k$-bonacci matrices with the help of $k$-bonacci matrix of order $k$ [see (\ref{3})] as follows.

\[
\begin{array}{ccc}
\sum_{n=0}^{k-1} F_{j+n,k}^{(k)} 
&=&\sum_{n=0}^{k-1} \left(f_{\lambda\mu+n,k}\right)
(\mbox{from equation}\;(\ref{3}))\\
&=&\left(\sum_{n=0}^{k-1} f_{\lambda\mu+n,k}\right)\\
&=&\left(f_{\lambda\mu+k,k}\right)\\
&=&F_{j+k,k}^{(k)}.
\end{array}
\]
\noindent Therefore we have the k-bonacci matrix of order $k$, we can construct the k-bonacci matrix of order $k^{2}$ 

\[
\begin{array}{ccc}
\sum_{n=0}^{k-1} F_{j+n,k}^{(k^{2})}
&=&\sum_{n=0}^{k-1}\left(F_{\lambda\mu+n,k}^{(k)}\right)
(\mbox{from}\;(\ref{3}))\\
&=&\left(\sum_{n=0}^{k-1}F_{\lambda\mu+n,k}^{(k)}\right)\\
&=&\left(F_{\lambda\mu+k,k}^{(k)}\right)\\
&=&F_{j+k,k}^{(k^{2})}.
\end{array}
\]

\noindent Similarly the k-bonacci matrix of order $k^{r}$ can be written as

\[
\begin{array}{ccc}
\sum_{n=0}^{k-1} F_{j+n,k}^{(k^{r})}
&=&\sum_{n=0}^{k-1} \left(F_{\lambda\mu+n,k}^{(k^{r-1})}\right)
(\mbox{from}\;(\ref{3}))\\
&=&\left(\sum_{n=0}^{k-1} F_{\lambda\mu+n,k}^{(k^{r-1})}\right)\\
&=&\left(F_{\lambda\mu+k,k}^{(k^{r-1})}\right)\\
&=&F_{j+k,k}^{(k^{r})}.
\end{array}
\]

\noindent More precisely, the k-bonacci matrix of order $k^{r}$ is given as
\[
F_{j+k,k}^{(k^{r})}=\sum_{n=0}^{k-1} F_{j+n,k}^{(k^{r})} \quad\text{for all j and r}, \quad\text {where}\quad j \geq 0, r \geq 1.
\]
%%%%%%%%%%%%%%%%%%%%%%%%%%%%%%%%%%%%%%%%%%%%%%%%%%%%%%%%%%%%%%%%%%%%%%%%%%%%%%%%%%%%%%%%%%%%%%%%%%%%%%%
\section{Properties of Higher Order Matrices}

Lemma $4.1$ can be obtained by using equation (\ref{one-e}) and Theorem $4$ of \cite{gab70}.

\begin{lemma}\label{fl}
\[
\sum_{j=0}^{n-1} f_{j,k}=\frac{1}{k-1}  
\left(f_{n+k-1,k}-f_{k-1,k}-\sum_{i=1}^{k-2} 
if_{n+k-2-i,k}\right).
\]
\end{lemma}

The {\bf Sum Formula} for the higher order $k$-bonacci matrices is given by the following Theorem.

\begin{theorem} \label{sum} 
If $n \geq 1,j \geq 0$ and $r \geq 1$, then

\[
\sum_{j=0}^{n-1} F_{j,k}^{(k^{r})}=\frac{1}{k-1}  
\left(F_{n+k-1,k}^{(k^{r})}-F_{k-1,k}^{(k^{r})}-
\sum_{i=1}^{k-2} iF_{n+k-2-i,k}^{(k^{r})}\right).
\]
\end{theorem}

\proof
We have the matrix from equation $(2.3)$
\[
\sum_{j=0}^{n-1} F_{j,k}^{(k^{r})}
=\sum_{j=0}^{n-1}\left(F_{j+k-\lambda-\mu+1,k}^{(k^{r-1})}\right)
\]
The summation on the right hand side of the above expression can be further simplified using Lemma 4.1
\[
\begin{array}{ccc}
&=&\frac{1}{k-1}\left(F_{n+k-1+k-\lambda-\mu+1,k}^{(k^{r-1})}-F_{k-1+k-\lambda-\mu+1,k}^{(k^{r-1})}-\sum_{i=1}^{k-2} iF_{n+k-2-i+k-\lambda-\mu+1,k}^{(k^{r-1})}\right)\\& & \\
&=&\frac{1}{k-1}\left(F_{n+2k-\lambda-\mu,k}^{(k^{r-1})}-F_{2k-\lambda-\mu,k}^{(k^{r-1})}-\sum_{i=1}^{k-2} iF_{n+2k-1-i-\lambda-\mu,k}^{(k^{r-1})}\right)\\& & \\
&=&\frac{1}{k-1}\left(F_{n+k-1,k}^{(k^{r})}-F_{k-1,k}^{(k^{r})}-\sum_{i=1}^{k-2}iF_{n+k-2-i,k}^{(k^{r})}\right).
\end{array}
\]
Thus, Theorem $4.2$ holds true for $n \geq 1,j \geq 0$ and $r \geq 1$.
\qed

Lemma $4.3$ is obtained by equation (\ref{one-e}) and Theorem $1$ of \cite{gab70}.

\begin{lemma}\label{sl}  
If $j \geq 2$, then
\[
f_{j+k-1,k}=2 f_{j+k-2,k}-f_{j-2,k}.
\]
\end{lemma}

\begin{theorem}
If $j \geq 2$, then
\[
F_{j+k-1,k}^{(k^{r})}=2 F_{j+k-2,k}^{(k^{r})}-F_{j-2,k}^{(k^{r})}. 
\]
\end{theorem}

\proof
Right hand side of Theorem $4.4$ can be written in the following form by equation $(2.3)$,
\[
\begin{array}{ccc}
2F_{j+k-2,k}^{(k^{r})}-F_{j-2,k}^{(k^{r})}
&=&2\left(F_{j+k-2+k-\lambda-\mu+1,k}^{(k^{r-1})}\right)-\left(F_{j-2+k-\lambda-\mu+1,k}^{(k^{r-1})}\right)\\
&=&2\left(F_{j+2k-1-\lambda-\mu,k}^{(k^{r-1})}\right)
-\left(F_{j-1+k-\lambda-\mu,k}^{(k^{r-1})}\right)
\end{array}
\]
In particular, using Lemma $4.3$, we obtain
\[
\begin{array}{ccc}
&=&\left(F_{j+2k-\lambda-\mu,k}^{(k^{r-1})}\right)\\
& & \\
&=&F_{j+k-1,k}^{(k^{r})}.
\end{array}
\]
\qed

\begin{lemma}\label{tl}
If $n \geq 1$, then
\begin{equation} \label{4}
\sum_{j=1}^{n} \frac{f_{j-1,k}}{2^{j}}
=1-\frac{f_{k+n,k}}{2^{n}}.
\end{equation}
\end{lemma}

\proof
(Proof by Mathematical Induction) In the initial step, we must verify that the above expression is true for $n=1$,
\[
\begin{array}{ccc}
\frac{f_{0,k}}{2}&=&\frac{f_{k,k}-\left(\sum_{t=1}^{k-1}f_{t,k}\right)}{2}\\
&=&\frac{2f_{k,k}}{2}-\frac{f_{k+1,k}}{2}\\
&=&1-\frac{f_{k+1,k}}{2}.
\end{array}
\]
\noindent The expression is clearly true. We assume that expression (\ref{4}) is true for $n=m$, i.e,
\[
\sum_{j=1}^{m}\frac{f_{j-1,k}}{2^{j}}=1-\frac{f_{k+m,k}}{2^{m}}
\]
\noindent In the next step, we will show that it holds for $n=m+1$, 
\[
\begin{array}{ccc}
\sum_{j=1}^{m+1}\frac{f_{j-1,k}}{2^{j}}
&=&1-\frac{f_{k+m},k}{2^{m}}+\frac{f_{m,k}}{2^{m+1}}\\
&=&\frac{2^{m+1}-f_{k+m,k}+f_{k+m,k}
-\left(\sum_{t=1}^{k}f_{m+t,k}\right)}{2^{m+1}}\\
&=&1-\frac{f_{k+m+1,k}}{2^{m+1}}.
\end{array}
\]
This is the required result.
\qed

\begin{theorem}
\[
\sum_{j=1}^{n} \frac{F_{j-1,k}^{(k^{r})}}{2^{j}}
=1-\frac{F_{k+n,k}^{(k^{r})}}{2^{n}}
\]
\end{theorem}

\proof
Theorem $4.6$ can be proved by using equation $(2.3)$ and Lemma $4.5$.
\[
\sum_{j=1}^{n} \frac{F_{j-1,k}^{(k^{r})}}{2^{j}}
=\sum_{j=1}^{n}\frac{\left(F_{j+k-\lambda-\mu,k}^{(k^{r-1})}\right)}{2^{j}}\\
=1-\frac{\left(F_{2k+n-\lambda-\mu+1,k}^{(k^{r-1})}\right)}{2^{n}}
\]
\qed

The following Lemma is obtained by equation (\ref{one-e}) and theorem $3(a)$ of \cite{gab70}.

\begin{lemma}\label{fol}
\[
\sum_{n=0}^{m} f_{kn+j+1,k}
=\sum_{n=-mk}^{k-1} f_{j-n,k}.
\]
\end{lemma}

\begin{theorem}\label{ft}
If $j \geq 0$ and $r \geq 1$, then
\[
\sum_{n=0}^{m} F_{kn+j+1,k}^{(k^{r})}
=\sum_{n=-mk}^{k-1} F_{j-n,k}^{(k^{r})}.
\]
\end{theorem}

\proof
From the matrix of $(2.3)$, we have
\[
\sum_{n=0}^{m} F_{kn+j+1,k}^{(k^{r})}
=\sum_{n=0}^{m} \left(F_{k(n+1)+j-\lambda-\mu+2,k}^{(k^{r-1})}\right)
\]
Which gives, using lemma $4.7$, the following.
\[
\begin{array}{ccc}
&=&\sum_{n=-mk}^{k-1} \left(F_{j-n+k-\lambda-\mu+1,k}^{(k^{r-1})}\right)\\
&=&\sum_{n=-mk}^{k-1} F_{j-n,k}^{(k^{r})}
\end{array}
\]
\qed

Lemma $4.9$ is obtained by equation (\ref{one-e}) and Theorem $3(b)$ of \cite{gab70}.

\begin{lemma}\label{fil}
If $m \geq 1$, then
\[
\sum_{n=1}^{m} f_{kn,k}
=\sum_{n=k(1-m)}^{k-1} f_{k-1-n,k}.
\]
\end{lemma}

\begin{theorem}\label{st}
If m is greater than 1, then
\[
\sum_{n=1}^{m} F_{kn,k}^{(k^{r})}
=\sum_{n=k(1-m)}^{k-1} F_{k-1-n,k}^{(k^{r})}
\]
\end{theorem}

\proof
With the help of equation $(2.3)$, the summation of the left hand side of Theorem $4.10$, can be written as 
\[
\sum_{n=1}^{m} F_{kn,k}^{(k^{r})}
=\sum_{n=1}^{m} \left(F_{k(n+1)-\lambda-\mu+1,k}^{(k^{r-1})}\right)
\]
Further, by Lemma $4.9$, we get
\[
\begin{array}{ccc}
&=\sum_{n=k(1-m)}^{k-1} \left(F_{2k-n-\lambda-\mu,k}^{(k^{r-1})}\right)\\
&=\sum_{n=k(1-m)}^{k-1} F_{k-1-n,k}^{(k^{r})}
\end{array}
\]
\qed

Lemma $4.11$ is obtained by equation (\ref{one-e}) and Theorem $3(c)$ of \cite{gab70}.
\begin{lemma}\label{sil}
\[
f_{km,k}-f_{0,k}
=\sum_{n=k(1-m)\;n \neq 0\;(mod\;k-1)}^{k-1} f_{k-1-n,k}.
\]
\end{lemma}

\begin{theorem}
If $m \geq 1$, then
\[
F_{km,k}^{(k^{r})}-F_{0,k}^{(k^{r})}=
\sum_{n=k(1-m)\;n \neq 0\; (mod\;k-1)}^{k-1} 
F_{k-1-n,k}^{(k^{r})}.
\]
\end{theorem}

\proof
The above theorem immediately follows from the Theorem \ref{st}. 
\qed

By Theorem $5$ of \cite{gab70} and by equation (\ref{one-e}), we can get the following Lemma.

\begin{lemma}\label{sel}
It states that if $n \geq 0$, then
\[
\sum_{j=0}^{n} f_{j,k}^{2}
=f_{n+1,k}f_{n,k}-\sum_{j=2}^{k-1} 
\sum_{i=0}^{n} f_{i,k} f_{i-j,k}.
\]
\end{lemma}

\begin{theorem}
If $n \geq 0$, then
\[
\sum_{j=0}^{n} (F_{j,k}^{(k^{r})})^{2}
=F_{n+1,k}^{(k^{r})}F_{n,k}^{(k^{r})}-
\sum_{j=2}^{k-1} \sum_{i=0}^{n} F_{i,k}^{(k^{r})}
F_{i-j,k}^{(k^{r})}.
\]
\end{theorem}

\proof
By equation $(2.3)$, right hand side of the theorem can be written as 
\[
=\left(F_{n+k-\lambda-\mu+2,k}^{(k^{r-1})}\right)
\left(F_{n+k-\lambda-\mu+1,k}^{(k^{r-1})}\right)
-\sum_{j=2}^{k-1} \sum_{i=0}^{n}
\left(F_{i-k-\lambda-\mu+1,k}^{(k^{r-1})}\right) 
\left(F_{i-j+k-\lambda-\mu+1,k}^{(k^{r-1})}\right)
\]
From Lemma $4.13$, we have 
\[
\begin{array}{cc}
&=\left(\sum_{j=0}^{n} F_{j+k-\lambda-\mu+1,k}^{(k^{r-1})}\right)^{2}\\ 
&=\sum_{j=0}^{n} (F_{j,k}^{(k^{r})})^{2}
\end{array}
\]
Thus, the theorem is proved.
\qed

\begin{theorem}\label{eig}
If $n \geq 1$, then
\begin{equation}\label{5}
f_{j+n,k}=2^{n}\sum_{i=1}^{k-n}f_{j-i,k}+(2^{n}-1)f_{j-k+(n-1),k}+(2^{n}-2)f_{j-k+(n-2),k}+(2^{n}-2^{2})f_{j-k+(n-3),k}+\ldots+2^{n-1}f_{j-k,k}.
\end{equation}
\end{theorem}

\proof
In the initial step, proof is carried out by induction on $n$. Clearly it is true for $n=1$, 
\[
\begin{array}{ccc}
&=&2\sum_{i=1}^{k-1} f_{j-i,k}+f_{j-k,k} \\
&=& 2 (f_{j-1,k}+\ldots+f_{j-k+1,k})+f_{j-k,k}\\
&=& f_{j+1,k}
\end{array}
\]
Thus we assume that the expression (\ref{5}) is true for $n=m$, so far it should hold true for $n=m+1$, that can be observe in the following expression
\[
\begin{array}{ccc}
f_{j+m+1,k}&=&2^{m+1}\sum_{i=1}^{k-(m+1)} f_{j-i,k}
+(2^{m+1}-1)f_{j-(k-(m+1)-1),k}+(2^{m+1}-2)f_{j-k+(m-1),k}+\ldots+2^{m}f_{j-k,k}\\
&=&2^{m+1}\sum_{i=1}^{k-m}f_{j-i,k}-f_{j-k+m,k}+2(2^{m}-1)f_{j-k+(m-1),k}+\ldots+2.2^{m-1}f_{j-k,k}\\
&=&2f_{j+m,k}-f_{j-k+m,k}
\end{array}
\]
Thus, by induction and the virtue of Lemma \ref{sl}, we have the desired result.
\qed

\begin{theorem} 
\begin{equation}\label{6}
F_{j+n,k}^{(k^{r})}=2^{n}\sum_{i=1}^{k-n} 
F_{j-i,k}^{(k^{r})}+(2^{n}-1)F_{j-k+(n-1),k}^{(k^{r})}
+(2^{n}-2)F_{j-k+(n-2),k}^{(k^{r})}
+\ldots+2^{n-1}F_{j-k,k}^{(k^{r})}.
\end{equation}
\end{theorem}

\proof
Proof of the Theorem $4.16$ is carried out by induction on $n$. For $n=1$, the expression (\ref{6}) takes the form
\[
=2\sum_{i=1}^{k-1} F_{j-i+k-\lambda-\mu+1,k}^{(k^{r-1})}+F_{j-\lambda-\mu+1,k}^{(k^{r-1})} 
\]
By equation (2.3), it takes the form
\[
\begin{array}{ccc}
&=& 2 (F_{j+k-\lambda-\mu,k}^{(k^{r-1})}+F_{j-1+k-\lambda-\mu,k}^{(k^{r-1})}+\ldots+F_{j-\lambda-\mu+2,k}^{(k^{r-1})})+F_{j-\lambda-\mu+1,k}^{(k^{r-1})} \\
&=& F_{j+k+2-\lambda-\mu,k}^{(k^{r-1})} \\
&=& F_{j+1,k}^{(k^{r})}
\end{array}
\]
Since result is true for $n=1$. We assume that expression (\ref{6}) is true for $n=m$. For $n=m+1$ it directly follows from the proof of Theorem \ref{eig}.
\qed

\begin{theorem}
If $r \geq 1$, then
\[
F_{j,k}^{(k^{r})}=F_{1,k}^{(k^{r})}
\left(\begin{array}{cccccc}
I_{k^{r-1}} & I_{k^{r-1}} & O_{k^{r-1}} & \cdots & O_{k^{r-1}} \\
I_{k^{r-1}} & O_{k^{r-1}} & I_{k^{r-1}} & \cdots & O_{k^{r-1}} \\
\vdots & \vdots &  \vdots & \vdots  & \vdots & \vdots \\
I_{k^{r-1}} & O_{k^{r-1}} & O_{k^{r-1}} & \cdots & I_{k^{r-1}} \\
I_{k^{r-1}} & O_{k^{r-1}} & O_{k^{r-1}} & \cdots & O_{k^{r-1}} \\
\end{array}
\right)^{j-1}
=F_{1,k}^{(k^{r})}(Q_{k}^{r})^{j-1}.
\]
\end{theorem}

\noindent Where $I_{k^{r-1}}$ is the identity matrix of order $k^{r-1}$
and $O_{k^{r-1}}$ is the null matrix of order $k^{r-1}$. Matrix $Q$ is the
combined representation of both matrices.

\proof
We have the matrix from $(2.3)$
\[
F_{j,k}^{(k^{r})}=F_{j+k-\lambda-\mu+1,k}^{(k^{r-1})}
\]
Above expression can be rewritten as
\[
\begin{array}{cc}
&=F_{j+k-\lambda-\mu,k}^{(k^{r-1})}(Q_{k}^{r})\\
&=F_{j+k-\lambda-\mu-1,k}^{(k^{r-1})}(Q_{k}^{r})^{2}\\
&=F_{k+2-\lambda-\mu,k}^{(k^{r-1})}(Q_{k}^{r})^{j-1}
\end{array}
\]
Hence, the theorem is proved.
\qed

%%%%%%%%%%%%%%%%%%%%%%%%%%%%%%%%%%%%%%%%%%%%%%%%%%%%%%%%%%%%%%%%%%%%%%%%%%%%%%%%%%%%%%%%%%%%%%%%%%%%%%%%%%
\section{Some Fibonacci Results of Matrices of Order $2^{r}$}

In this section, we obtain some new Fibonacci identities for matrices of order $2^{r}$. In order to get these identities, we have used Fibonacci \cite{che01} and Lucas matrices \cite{chen01}.

For the proof of Lemmas $5.1, 5.4, 5.6\quad\text{and}\quad5.8$ see the paper \cite{Vaj89}. In particular, Lemma $5.1$ contains the well-defined identities of Fibonacci and Lucas numbers.

\begin{lemma}
(i)$F_{m+n} = F_{m-1}F_{n}+F_{m}F_{n+1}$
(ii) $F_{m+n}+F_{m+n+2} = L_{m+n+1}$
(iii) $L_{n+1}+L_{n+3} = 5F_{n+2}$\\
for all integers n
\end{lemma}

\begin{theorem}\label{11}
\[
(i) L_{m+n}^{(2^{r})}+L_{m+n+2}^{(2^{r})}=5F_{m+n+1}^{(2^{r})}
\]
\[
(ii) F_{m+n}^{(2^{r})}+F_{m+n+2}^{(2^{r})}=L_{n+m+1}^{(2^{r})}
\]
\end{theorem}

\proof
First identity of the Theorem 5.2 can be proved by applying induction on r, therefor for $r=1$, we have   
\[
L_{n+m}^{(2)}+L_{n+m+2}^{(2)}=
\left(\begin{array}{cc} L_{n+m+1} & L_{n+m}\\
L_{n+m} & F_{n+m-1}
\end{array}\right)
+ 
\left(\begin{array}{cc} L_{n+m+3} & L_{n+m+2}\\
L_{n+m+2} & L_{n+m+1}
\end{array}\right)
=\left(\begin{array}{cc} 5F_{n+m+2} & 5F_{n+m+1}\\
5F_{n+m+1} & 5F_{n+m} 
\end{array}\right)
=5F_{n+m+1}^{(2)} 
\]
\noindent Which is true. Thus, we assume that it is true for $r=p$, i.e.
$$ L_{n+m}^{(2^{p})} +  L_{n+m+2}^{(2^{p})}
= 5F_{n+m+1}^{(2^{p})}$$
Therefore, we should prove it for $r=p+1$,
\[
L_{n+m}^{(2^{p+1})}+L_{n+m+2}^{(2^{p+1})}
=\left(\begin{array}{cc} L_{n+m+1}^{(2^{p})} & L_{n+m}^{(2^{p})}\\
L_{n+m}^{(2^{p})} & L_{n+m-1}^{(2^{p})}
\end{array}\right)
+ 
\left(\begin{array}{cc} L_{n+m+3}^{(2^{p})} & L_{n+m+2}^{(2^{p})}\\
L_{n+m+2}^{(2^{p})} & L_{n+m+1}^{(2^{p})}
\end{array}\right)
\]
\[
=\left(\begin{array}{cc} 5F_{n+m+2}^{(2^{p})} & 5F_{n+m+1}^{(2^{p})}\\
5F_{n+m+1}^{(2^{p})} & 5F_{n+m}^{(2^{p})} 
\end{array}\right)
\]
\[
=5F_{n+m+1}^{(2^{p+1})} 
\]
Hence, we have the result.
\qed

(ii) Second identity of the theorem can be proved in a similar way.

\begin{theorem}\label{12}
If $r=2t$ and $t \geq 1$, then
\[
(i) 5^{\frac{r}{2}}F_{m+n}^{(2^{r})}=F_{m-1}^{(2^{r})}F_{n}^{(2^{r})}
+F_{m}^{(2^{r})}F_{n+1}^{(2^{r})}
\]
If $r=2s+1$ and $s \geq 0$, then
\[
(ii) 5^{\frac{r-1}{2}}L_{m+n}^{(2^{r})}=F_{m-1}^{(2^{r})}F_{n}^{(2^{r})}
+F_{m}^{(2^{r})}F_{n+1}^{(2^{r})}
\]
\end{theorem}
 
\proof
Even $r$ can be obtained from odd $r$ by applying induction on t. On taking $t=1$, r becomes $r=2$ in the identity (i) of the Theorem 5.3, therefore we have 
\[
F_{m-1}^{(2^{2})}F_{n}^{(2^{2})}+F_{m}^{(2^{2})}F_{n+1}^{(2^{2})}
\]
\[
= \left(\begin{array}{cc} 
F_{m}^{(2)} & F_{m-1}^{(2)}\\
F_{m-1}^{(2)} & F_{m-2}^{(2)}
\end{array}\right)
\left(\begin{array}{cc} 
F_{n+1}^{(2)} & F_{n}^{(2)}\\
F_{n}^{(2)} & F_{n-1}^{(2)}
\end{array}\right)
+
\left(\begin{array}{cc} 
F_{m+1}^{(2)} & F_{m}^{(2)}\\
F_{m}^{(2)} & F_{m-1}^{(2)}
\end{array}\right)
\left(\begin{array}{cc} 
F_{n+2}^{(2)} & F_{n+1}^{(2)}\\
F_{n+1}^{(2)} & F_{n}^{(2)}
\end{array}\right)
\]
\noindent By multiplying and simplifying these matrices, we get (see Theorem 5.2)
\[ 
=\left(\begin{array}{cc} 
L_{n+m}^{(2)} & L_{n+m-1}^{(2)}\\
L_{n+m-1}^{(2)} & L_{n+m-2}^{(2)}
\end{array}\right)
+
\left(\begin{array}{cc} 
L_{n+m+2}^{(2)} & L_{n+m+1}^{(2)}\\
L_{n+m+1}^{(2)} & L_{n+m}^{(2)}
\end{array}\right)
\]
\[
=5F_{n+m}^{(2^{2})}
\]
\noindent Let suppose that Theorem $5.3(i)$ is true for $t=p$. We need to show it true for $t=p+1$,
\[
F_{m-1}^{(2^{2p+2})}F_{n}^{(2^{2p+2})}
+F_{m}^{(2^{2p+2})}F_{n+1}^{(2^{2p+2})}
\]
\[
= \left(\begin{array}{cc} 
F_{m}^{(2^{2p+1})} & F_{m-1}^{(2^{2p+1})}\\
F_{m-1}^{(2^{2p+1})} & F_{m-2}^{(2^{2p+1})}
\end{array}\right)
\left(\begin{array}{cc} 
F_{n+1}^{(2^{2p+1})} & F_{n}^{(2^{2p+1})}\\
F_{n}^{(2^{2p+1})} & F_{n-1}^{(2^{2p+1})}
\end{array}\right)
+\left(\begin{array}{cc} 
F_{m+1}^{(2^{2p+1})} & F_{m}^{(2^{2p+1})}\\
F_{m}^{(2^{2p+1})} & F_{m-1}^{(2^{2p+1})}
\end{array}\right)
\left(\begin{array}{cc} 
F_{n+2}^{(2^{2p+1})} & F_{n+1}^{(2^{2p+1})}\\
F_{n+1}^{(2^{2p+1})} & F_{n}^{(2^{2p+1})}
\end{array}\right)
\]
\[
= 5^{p}\left(\begin{array}{cc} 
L_{n+m}^{(2^{2p+1})} & L_{n+m-1}^{(2^{2p+1})}\\
L_{n+m-1}^{(2^{2p+1})} & L_{n+m-2}^{(2^{2p+1})}
\end{array}\right)
+ 
5^{p} \left(\begin{array}{cc} 
L_{n+m+2}^{(2^{2p+1})} & L_{n+m+1}^{(2^{2p+1})}\\
L_{n+m+1}^{(2^{2p+1})} & L_{n+m}^{(2^{2p+1})}
\end{array}\right)
\] 

\[
=5^{(p+1)}F_{n+m}^{(2^{2p+2})}
\]
We showed that it is true $t=p+1$.
\qed

Identity(ii) of Theorem 5.3 can be proved by induction on $s$, therefore when $s=0$, r becomes $r=1$, identity can be written as
\[
F_{m-1}^{(2)}F_{n}^{(2)}+F_{m}^{(2)}F_{n+1}^{(2)}
\]
\noindent Then, we have  

\[
\left(\begin{array}{cc} 
F_{n+m} & F_{n+m-1}\\
F_{n+m-1} & F_{n+m-2}
\end{array}\right)
+ 
\left(\begin{array}{cc} 
F_{n+m+2} & F_{n+m+1}\\
F_{n+m+1} & F_{n+m}
\end{array}\right)
\]
\noindent which becomes by lemma $5.1$
\[
\left(\begin{array}{cc} 
L_{n+m+1} & L_{n+m}\\
L_{n+m} & L_{n+m-1}
\end{array}\right)
=L_{n+m}^{(2)}
\]
\noindent Which is absolutely true. For $s=1$, it is also true 
\[
F_{m-1}^{(2^{3})}F_{n}^{(2^{3})}+F_{m}^{(2^{3})}F_{n+1}^{(2^{3})}
\]
\[
\begin{array}{cc}
=&5(F_{m+n-1}^{(2^{2})}+F_{m+n+1}^{(2^{2})})\\
=&5L_{n+m}^{(2^{3})}
\end{array}
\]
\noindent Therefore we suppose that Theorem (12(ii)) is true for $s=q-1$. We need to show that it is also true for $s=q$
\[
F_{m-1}^{(2^{2q+1})}F_{n}^{(2^{2q+1})}
+F_{m}^{(2^{2q+1})}F_{n+1}^{(2^{2q+1})}
\]
\[
F_{m-1}^{(2^{2q})}F_{m-1}^{(2)}
F_{n}^{(2^{2q})}F_{n}^{(2)}
+F_{m}^{(2^{2q})}F_{m}^{(2)}
F_{n+1}^{(2^{2q})}F_{n+1}^{(2)}
\]
\noindent After simplifying, we have the desire result.
\[
=(F_{n+m-1}^{(2)}5^{q}F_{n+m-1}^{(2^{2q})}
+F_{n+m+1}^{(2)}5^{q}F_{n+m+1}^{(2^{2q})})
\]
\[
=5^{q}L_{n+m}^{(2^{2q+1})}
\]
\qed

\begin{lemma}
\[ 
F_{n}+L_{n}=2F_{n+1}
\]
\end{lemma}

\begin{theorem}
\[
F_{n}^{(2^{r})}+L_{n}^{(2^{r})}=2F_{n+1}^{(2^{r})}
\]
\end{theorem}

\proof
(Proof by Mathematical Induction on $r$) For $r=1$, we have
\[
F_{n}^{(2)}+L_{n}^{(2)}=
\left(\begin{array}{cc} F_{n+1} & F_{n}\\
F_{n} & F_{n-1}
\end{array}\right)
+ 
\left(\begin{array}{cc} L_{n+1} & L_{n}\\
L_{n} & L_{n-1}
\end{array}\right)
=\left(\begin{array}{cc} 2F_{n+2} & 2F_{n+1}\\
2F_{n+1} & 2F_{n} 
\end{array}\right)
=2F_{n+1}^{(2)} 
\]
\noindent Which is clearly true. Assume that the identity is true for $r=n$, i.e.
\begin{equation}
F_{n}^{(2^{r})}+L_{n}^{(2^{r})}=2F_{n+1}^{(2^{r})}
\end{equation}
\noindent Therefore, for $r=n+1$, we need to show that it is true
\[
F_{n}^{(2^{n+1})}+L_{n}^{(2^{n+1})}
=\left(\begin{array}{cc} F_{n+1}^{(2^{n})} & F_{n}^{(2^{n})}\\
F_{n}^{(2^{n})} & F_{n-1}^{(2^{n})}
\end{array}\right)
+ 
\left(\begin{array}{cc} L_{n+1}^{(2^{n})} & L_{n}^{(2^{n})}\\
L_{n}^{(2^{n})} & L_{n-1}^{(2^{n})}
\end{array}\right)
\]
\[
=\left(\begin{array}{cc} 2F_{n+2}^{(2^{n})} & 2F_{n+1}^{(2^{n})}\\
2F_{n+1}^{(2^{n})} & 2F_{n}^{(2^{n})} 
\end{array}\right)
\] 
\[
=2F_{n+1}^{(2^{n+1})} 
\]
Hence, we have the result.
\qed

\begin{lemma}
\[ 
F_{n+1}^{2}+F_{n}^{2}=F_{2n+1}
\]
\end{lemma}

\begin{theorem}
If $r=2t$ and $t \geq 1$, then
\[
(i)  5^{\frac{r}{2}}F_{2n+1}^{(2^{r})}
=(F_{n+1}^{(2^{r})})^{2}+(F_{n}^{(2^{r})})^{2}
\]
If $r=2s+1$ and $s \geq 0$, then
\[
(ii)  5^{\frac{r-1}{2}}L_{2n+1}^{(2^{r})}
=(F_{n+1}^{(2^{r})})^{2}+(F_{n}^{(2^{r})})^{2}
\]
\end{theorem}

\proof
It is similar to the proof of Theorem \ref{12}.
\qed

\begin{lemma}
\[ 
F_{n+1}^{2}-F_{n}^{2}=F_{n+2}F_{n-1}
\]
\end{lemma}

\begin{theorem}
If $r=2t$ and $t \geq 1$, then
\[
(i) 5^{\frac{r}{2}}L_{2n+1}^{(2^{r})}
=(F_{n+1}^{(2^{r})})^{2}-(F_{n}^{(2^{r})})^{2}
\]
If $r=2s+1$ and $s \geq 0$, then
\[
(ii)  5^{\frac{r-1}{2}}F_{2n+1}^{(2^{r})}
=(F_{n+1}^{(2^{r})})^{2}-(F_{n}^{(2^{r})})^{2}
\]
\end{theorem}

\proof
It follows easily from the proof of Theorem \ref{12}.
\qed

\begin{lemma}
\[
\sum_{i=1}^{n}F_{i}^{2}=F_{n}F_{n+1}
\]
\end{lemma}

\begin{theorem}
If $r=2s+1$, $s \geq 0$ and $n \geq 1$, then
\[
(i)  \sum_{i=1}^{n} 
\left(F_{i}^{(2^{r})}\right)^{2}
=5^{\frac{r-1}{2}}\left(F_{2n+1}^{(2^{r})}
-F_{1}^{(2^{r})}\right) 
\]
If $r=2t$, $t \geq 1$ and $n \geq 1$, then
\[
(ii)  \sum_{i=1}^{n} 
\left(F_{i}^{(2^{r})}\right)^{2}
=5^{\frac{r-2}{2}}\left(L_{2n+1}^{(2^{r})}
-L_{1}^{(2^{r})}\right) 
\]
\end{theorem}

\proof
Since
\[
\left(F_{i}^{(2)}\right)^{2}
= \left(\begin{array}{cc} 
F_{i+1}^{2}+F_{i}^{2} & F_{i}(F_{i+1}+F_{i-1})\\
F_{i}(F_{i+1}+f_{i-1}) & F_{i}^{2}+F_{i-1}
\end{array}\right)
\]
\[
=\left(\begin{array}{cc} 
F_{2i+1} & F_{i}L_{i}\\
F_{i}L_{i} & F_{2i-1}
\end{array}\right)
\]
\noindent By the formula $13$ on the page $177$ of \cite{Vaj89}, we get 
\[
=\left(\begin{array}{cc} 
F_{2i+1} & F_{2i}\\
F_{2i} & F_{2i-1}
\end{array}\right)
\]
\noindent By mathematical induction on $s$, by taking $s=0$, r becomes $r=1$, we have
\[
\sum_{i=1}^{n}\left(F_{i}^{(2)}\right)^{2}
\]
\[
=\sum_{i=1}^{n}\left(\begin{array}{cc} 
F_{2i+1} & F_{2i}\\
F_{2i} & F_{2i-1}
\end{array}\right)
\]
\[
=\sum_{i=1}^{n}F_{2i}^{(2)}
=\sum_{i=2}^{2n}F_{i}^{(2)}
=F_{2n+1}^{(2)}-F_{1}^{(2)}
\]
\noindent It is directily followed by Theorem $5$ of \cite{che01}.\\
\noindent Since above expression is true for $s=0$. We assume that it is true for $s=p-1$. We need to prove it for $s=p$,
\[
\sum_{i=1}^{n}\left(F_{i}^{(2^{2p+1})}\right)^{2}
\]
\[
=\sum_{i=1}^{n}\left(\begin{array}{cc} 
F_{i}^{(2^{2p})} & F_{i}^{(2)}
\end{array}\right)^{2}
\]
\[
=\sum_{i=1}^{n}5^{p}F_{2i}^{(2^{2p})}F_{2i}^{(2)}
=5^{p}\sum_{i=1}^{n}F_{2i}^{(2^{2p+1})}
=5^{p}\sum_{i=2}^{2n}F_{i}^{(2^{2p+1})}
\]
\[
=5^{p}\left(F_{2n+1}^{(2^{2p+1})}-F_{1}^{(2^{2p+1})}\right)
\]
We have proved the first identity of Theorem 5.11.

(ii) Second identity of the Theorem can be prove in a similar way.

\qed

%%%%%%%%%%%%%%%%%%%%%%%%%%%%%%%%%%%%%%%%%%%%%%%%%%%%%%%References%%%%%%%%%%%%%%%%%%%%%%%%%%%%%%%%

\medskip

\noindent MSC2010: 11B39, 11C20, 15A24, 15A36


\begin{thebibliography}{99}

\bibitem{che01}
Z. Cheng-Heng,
{\it On the Fibonacci matrix of order $2^k$ sequence $\left\{F_n^{(2^k)}\right\}$},
{Ars. Combin.},
{\bf 59} (2001), 45--54.

\bibitem{chen01}
Z. Cheng-Heng,
{\it On the Generalized Fibonacci matrix of order $2^k$ sequence $\left\{G_n^{(2^k)}\right\}$},
{Ars. Combin.},
{\bf 59} (2001), 215--224.

\bibitem{Daz96}
L. Dazheng,
{\it Fibonacci matrices},
{The Fibonacci Quart.},
{\bf 37} (1999), 14--20.

\bibitem{gab70}
H. Gabai,
{\it Generalized Fibonacci $k$-sequences},
{The Fibonacci Quart.},
{\bf 8} (1970), no. 1, 31--38.

\bibitem{gou81}
H. W. Gould,
{\it A history of the Fibonacci Q-matrix and a higher-dimensional problem},
{The Fibonacci Quart.},
1981.

\bibitem{gra97}
F. Graham,
{\it The Singularity of Fibonacci Matrices},
{The Mathematical Gazette},
{\bf 81} (1979), no. 491, 295--298.

\bibitem{Hay99}
B. Hayes,
{\it The Vibonacci numbers},
{American Scientist}, 
{\bf 87} (1999), no. 4, 296--301.

\bibitem{hora61}
A. F. Horadam 
{\it A Generalized Fibonacci sequence}, 
{ Amer. Math. Mthly},
{\bf 68} (1961), 455--459.

\bibitem{Ivi}
J. Ivie,
{\it A general $Q$-Matrix}, 
{Fibnoacci Quart.},
{\bf 3} (1972), no. 10, 255--261.

\bibitem{Kil}
E. Kilic,
{\it The generalized order - k Fibonacci – Pell sequence by
matrix methods},
{Journal of Computational and Applied Mathematics},
{\bf 209} (2007) 133 –- 145.

\bibitem{Lee00}
G.-Y. Lee,
{\it k-Lucas numbers and associated bipartite graphs}, 
{ Linear Algebra and its applications},
{\bf 320} (2000), 51--61.

\bibitem{Lee97}
G.-Y. Lee, S.-G. Lee and H.-G. Shin,
{\it On the k-Generalized Fibonacci matrix $Q_{k}$}, 
{ Linear Algebra and its applications},
{\bf 251} (1997), 73--88.

\bibitem{Mei}
A. M. Meinke,
{\it FIBONACCI NUMBERS AND ASSOCIATED MATRICES},
{A thesis submitted to Kent State University},
(August 2011).

\bibitem{Mil60}
E. P. Miles,
{\it Generalized Fibonacci numbers and associated matrices},
{ Amer. Math. Monthly},
{\bf} (1960), 745--752.

\bibitem{nnm02}
N. N. Vorobʹev, N. N. Vorobiev and M. Martin,
{\it Fibonacci numbers},
{Birkhäuser, 2002}.

\bibitem{RaS}
J. L. Ram´ırez and V. F. Sirvent,
{\it A Generalization of the k-bonacci Sequence
from Riordan Arrays},
{ The electronic journal of combinatorics}, 
{\bf 22(1)} (2015), 1--38.

\bibitem{Sch97}
M. R. Schroeder, 
{\bf Number Theory in science and communication}, 
Springer, 1997.

\bibitem{Sir97}
V. F. Sirvent,
{\it A semigroup associated with the k-bonacci numbers with 
dynamic interpretation},
{The Fibonacci Quart.},
{\bf 35} (1997), no. 4, 335--340. 

\bibitem{Wh71}
J. E. Walton and A. F. Horadam,
{\it Some properties of certain generalized Fibonacci matrices}, 
{\bf 9.3} (1971), 264--276.

\bibitem{Vaj89}
S. Vajda,
{\bf Fibonacci $\&$ Lucas numbers, and the Golden section}, 
Ellis Horwood Limited, 1989.
\end{thebibliography}
\end{document}